\documentclass[a4paper,11pt]{amsart}
\usepackage{amssymb, graphics}
\input epsf
\newtheorem{lemma}{Lemma}

\newtheorem{thm}{Theorem}
\newtheorem{cor}{Corollary}
\theoremstyle{definition}
\newtheorem{rem}{Remark}

\newtheorem{conj}{Conjecture}
\newcounter{numl}

\newcommand{\labelnuml}{\textup{(\roman{numl})}}
\newenvironment{numlist}{\begin{list}{\labelnuml}%
{\usecounter{numl}\setlength{\leftmargin}{0pt}%
\setlength{\itemindent}{2\parindent}%
\setlength{\itemsep}{\smallskipamount}\def
\makelabel ##1{\hss \llap {\upshape ##1}}}}{\end{list}}

\newcommand{\C}{{\mathbb C}}

\newcommand{\cO}{{\mathcal O}}

\newcommand{\symprod}{\mathbin{\raise1pt\hbox{$\scriptstyle\bigcirc$}}}

\begin{document}
\title[K\"ahler metrics of constant scalar curvature on ruled 
surfaces]
{A remark on K\"ahler metrics of constant scalar curvature on ruled 
complex surfaces}
\author[V. Apostolov]{Vestislav Apostolov}
\address{Vestislav Apostolov \\ D{\'e}partement de Math{\'e}matiques\\
UQAM\\ C.P. 8888 \\ Succ. Centre-ville \\ Montr{\'e}al (Qu{\'e}bec) \\
H3C 3P8 \\ Canada}
\email{apostolo@math.uqam.ca}
\author[C. T\o nnesen-Friedman]{Christina W.~T\o nnesen-Friedman}
\address{Christina W. T\o nnesen-Friedman\\
Department of Mathematics\\ Union College\\
Schenectady\\  New York  12308\\ USA }
\email{tonnesec@union.edu}
\thanks{The first
author was supported in part by NSERC grant
OGP0023879, and the second author by the Union College Faculty Research Fund.}

\begin{abstract} In this note we point out how some recent 
developments in the theory of constant scalar curvature K\"ahler 
metrics can be used to clarify the existence issue for such metrics 
in the special case of (geometrically) ruled complex surfaces. 
\end{abstract}
\maketitle

\begin{center}
\emph{Dedicated to Professor Paul Gauduchon on his 60th birthday.}
\end{center}

\section{Introduction}
Let $(M,J)$ be a compact complex $m$-dimensional manifold and $\Omega \in 
H^{1,1}_{\mathbb R}(M) = H_{dR}^2(M, {\mathbb R}) \cap H^{1,1}(M, {\mathbb 
C})$ be a K\"ahler class on $M$.  We say that a K\"ahler metric $g$ 
belongs to $\Omega$ if the corresponding fundamental $(1,1)$-form 
$\omega(\cdot , \cdot) = g(J\cdot, \cdot)$ represents $\Omega$. 
Following Calabi~\cite{cal-one}, consider the  functional ${\rm Cal}(g) 
= \int_M s_g^2 d\mu_g$, where $s_g$ denotes the scalar 
curvature of $g$ and $d\mu_g= (1/m!) \omega^m$ is the volume form of 
$g$. Then $g$ is  called {\it extremal} if it is a critical 
point of ${\rm Cal}$ acting on the set of all K\"ahler metrics in 
$\Omega$. The corresponding Euler-Lagrange equation shows that $g$ is extremal
if and only if $J {\rm grad}_g  s_g$ is a Killing vector field on $(M,g)$. In 
particular, K\"ahler-Einstein metrics and  constant 
scalar curvature (CSC) K\"ahler metrics are extremal.

Obstructions to the existence of extremal or CSC metrics have been 
known since the pioneering work of Calabi~\cite{cal-two}. Indeed, Calabi 
showed that an extremal  K\"ahler metric is necessarily 
invariant under the action of a maximal compact subgroup 
of the connected component ${\rm Aut}_0(M,J)$ of the automorphism 
group ${\rm Aut} (M,J)$ of $(M,J)$. Furthermore, according to 
Futaki~\cite{futaki}, the scalar curvature of an extremal K\"ahler 
metric is constant if and only if the Futaki invariant identically 
vanishes at $\Omega$. The Matsushima--Lichnerowicz theorem~\cite{lichne,matsushima}, on the 
other hand, predicts that ${\rm Aut}_0(M,J)$ must be reductive if a 
CSC K\"ahler metric exists. Combining all these facts together, one 
easily finds lots of examples of K\"ahler manifolds which do not 
admit extremal K\"ahler metrics at all. However, besides these 
classical obstructions very little was known until recently for the 
existence and the uniqueness of extremal K\"ahler metrics  on a given 
compact complex manifold. 

Substantial progress was recently made after deep works of 
Donaldson~\cite{Do-one}, Chen--Tian~\cite{CT}, Luo~\cite{luo}, 
Mabuchi~\cite{mab-one,mab-two}, Paul--Tian~\cite{PT}, 
Zhang~\cite{zhang},   which relate the existence and the uniqueness 
issues for extremal K\"ahler metrics in Hodge K\"ahler classes with various notions of stability 
of polarized projective varieties. Adopting the terminology of 
\cite{RT}, the precise conjecture concerning CSC K\"ahler metrics is 
that K-polystability for a polarized (Hodge) K\"ahler 
manifold $(M,J, \Omega)$ should be equivalent to the existence of a 
CSC K\"ahler metric in the K\"ahler class $\Omega$. In one direction, namely that the existence of CSC K\"ahler metric implies K-stability, this conjecture has almost been proved.  
The only caveat is the apparent difference in definitions of 
stability used by different authors; in the notation of \cite{RT}, 
the results in \cite{CT,Do-one,mab-one}  imply  that if 
a polarized projective variety $(M,J, \Omega)$ admits a CSC 
K\"ahler metric in $\Omega$, then it must be K-semistable (where the 
conjectured K-polystability is a yet stronger condition).

While one immediate consequence of these recent 
works is the realization that  existence of extremal or CSC 
K\"ahler metrics depends in an essential way on the choice of  a 
K\"ahler 
class (or polarization) $\Omega$ (see \cite{Do-two,RT,christina-one} for 
examples),  the most exciting application  is the proof of the uniqueness:
\begin{thm}{\rm \cite{CT,Do-one,mab-one}}\label{uniqueness} Let $(M,J)$ be a compact  
complex manifold with a fixed K\"ahler class $\Omega$. Then any two 
extremal K\"ahler metrics in $\Omega$ must be biholomorphically 
isometric.
\end{thm}
This important result has been first established in \cite{Do-one}, under the assumptions that the K\"ahler class $\Omega$ is rational and the automorphism group ${\rm Aut}(M,J)$ is discrete; the latter assumption has been dropped in \cite{mab-one}, while in its most general version the above theorem was proved in~\cite{CT}.

The existence issue for CSC metrics, on the other hand, remains 
mysterious, mainly because it is difficult to verify whether a given 
polarized K\"ahler manifold satisfies the different notions of 
stability. A notable progress in this direction was recently made 
by Ross--Thomas \cite{RT} who introduced the more tractible 
notion of {\it slope} stability of a polarized projective 
variety. They show  that K-(semi)stability implies slope (semi)stability; thus the
existence of CSC metric in a given polarization is 
related to the relevant notion of slope polystability.

One special case where the numerical criterion of Ross--Thomas can be 
effectively applied is the class of (geometrically) ruled complex 
manifolds, which by definition are projectivisations $P(E)$ of 
holomorphic vector bundles over polarized smooth complex manifolds. 
In this situation, stability of the projective manifold $P(E)$ is 
related to the classical (slope) stability of the underlying 
holomorphic bundle $E$.  

Recall that the {\it slope} of a holomorphic vector bundle $E$ of rank $r$ over a compact K\"ahler manifold $(S,\omega_S)$ of complex dimension $n$ is the number $$\mu(E) = \frac{1}{r} \int_S c_1(E)\wedge \omega_S^{n-1}.$$
The vector bundle $E$ is called {\it stable} (resp.~{\it semistable}) if $\mu(F) < \mu(E)$ (resp.~$\mu(F) \le \mu(E)$) for any proper coherent sub-sheaf $F\subset E$. The bundle $E$ is {\it polystable} if it decomposes as a direct sum of stable vector bundles with the same slope. It is a well-known fact (see e.g. \cite{kobayashi}) that for vector bundles `stability' $\Rightarrow$ `polystability' $\Rightarrow$ `semistability'.  Note also that any holomorphic line bundle is stable.

At the current stage, the most general 
result  from \cite{RT}, combined with the theory of CSC K\"ahler metrics alluded to above,  show that  {\it if $E$ is not semistable, then $P(E)$ admits 
K\"ahler classes with no CSC K\"ahler metrics}.  

A partial converse is proved by Hong~\cite{hong}, who shows that if $E$ is a 
stable vector bundle over a Hodge K\"ahler manifold of CSC, such that 
$P(E)$ has discrete automorphism group, then 
suitable K\"ahler classes on $P(E)$ do admit CSC K\"ahler metrics.  

In between these results one also finds many {\it explicit} 
examples~\cite{ACGCT,guan,hwang,hwang-singer,koi-sak-one,koi-sak-two} of 
CSC K\"ahler  metrics on certain projective bundles. The most general assertion~\cite{ACGCT} is that the ruled manifold $M= P(E_1\oplus E_2) \to S$, where $S= \prod_{a=1}^N S_a$ is a product of K\"ahler CSC Hodge manifolds $(S_a, \omega_a)$, and $E_1$ and $E_2$ are projectively-flat (and therefore polystable~\cite{kobayashi}) Hermitian vector bundles of ranks $r_1$ and $r_2$, such that $c_1(E_2)/r_2 - c_1(E_1)/r_1 = \sum_{a=1}^N \varepsilon_a [\omega_a/2\pi]$  for $\varepsilon_a = \pm 1$, admits a CSC K\"ahler metric provided that not all of the $\varepsilon_a$'s have the same sign. This latter condition is equivalent to the assumption that the underlying vector bundle $E=E_1 \oplus E_2$ is polystable with respect to a suitable K\"ahler class on $S$. A notable feature of all the (non) existence results is that they depend, in an essential way,  on the choice of a K\"ahler class (or polarization) 
on both the base and the total space $P(E)$.

The situation is nicer if we consider a geometrically ruled complex 
manifold $P(E)$ over a compact complex curve. In this case, 
stability of the holomorphic vector bundle $E$ is manifestly independent of 
the choice of a K\"ahler metric on the base curve, and one could 
speculate that the notion of K-stability of the projective manifold $P(E)$ 
should be independent of the specific polarization. Indeed,  it follows from~\cite{RT} that if the variety $P(E)$ is K-semistable with respect to some 
polarization, then $E$ is a semistable vector bundle. As far 
as the existence of CSC K\"ahler metrics on such manifolds is 
concerned, one can go even further and observe that, 
when the manifold $P(E)$ has no non-trivial holomorphic vector fields with zeroes, 
the integrality condition on the 
K\"ahler class (i.e. the assumption that it defines a polarization) 
is inessential: The LeBrun-Simanca deformation 
theorem~\cite[Theorem~7]{Le-Sim} then tells us that if a K\"ahler class 
$\Omega$ admits a CSC K\"ahler metric, then, by denseness, 
a nearby {\it rational} K\"ahler class (and thus, by rescaling, a Hodge K\"ahler 
class) will as well. In particular, the general existence conjecture for 
CSC K\"ahler metrics alluded to above becomes, in the case of 
geometrically ruled manifolds over curves,  the following:
\begin{conj} {\it Let $(M,J)=P(E)$ be a geometrically ruled complex 
manifold over a compact complex curve.  Then $(M,J)$ admits a CSC 
K\"ahler metric in some {\rm (}and hence any{\rm )} K\"ahler class if and only if $E$ is polystable.} 
\end{conj}

\begin{rem} \noindent   (1) The `if' part is a well-known consequence 
of the classical 
theory in \cite{NS}. Indeed,  if $E$ is a polystable bundle over a 
curve, then $M=P(E)$ admits a locally-symmetric  CSC K\"ahler metric 
in each K\"ahler class. By Theorem~\ref{uniqueness}, these are the 
only CSC K\"ahler metrics on $M$. Thus Conjecture~1 can be 
reformulated as follows:~{\it any 
CSC K\"ahler metric on a geometrically ruled complex manifold over a 
compact complex curve is necessarily  locally-symmetric} (see also 
\cite[Lemma 8]{fujiki} and \cite{LeB}).

(2) In the case when $M=P(E)$ is a geometrically ruled complex manifold over ${\mathbb 
C} P^1$, $E$ splits as a direct sum of line bundles and the conjecture  follows from the 
Matsushima--Lichnerowicz theorem,  see e.g. \cite[Proposition~3]{ACGCT}.

(3) The main missing piece, to prove the above conjecture ---at least for rational K\"ahler classes--- by combining 
the results of \cite{RT} with the general theory \cite{CT,Do-one,mab-two}, is the fact that at the 
current stage it has only been shown that the existence of a CSC K\"ahler metric 
implies K-semistability (but not K-polystability).
\end{rem}

In this Note, we want to explain how 
Theorem~\ref{uniqueness} above, combined with  the  Matsushima--Lichnerowicz theorem~\cite{lichne,matsushima}, the Futaki obstruction theory~\cite{futaki}, and a deformation argument from 
\cite{fujiki,Le-Sim},  proves Conjecture~1 in the case when the 
rank of $E$ is two.

\begin{thm}\label{main} Let $(M,J)$ be a geometrically ruled complex 
surface, i.e. $M=P(E)$ is  the projectivisation of a holomorphic rank 
2 vector bundle $E$ over a compact complex curve $\Sigma_{\rm {\bf 
g}}$ of genus ${\rm {\bf g}}$. Then $M$ admits a CSC K\"ahler metric 
if and only if the holomorphic vector bundle $E$ is polystable {\rm (}i.e. 
is stable or is the direct sum of two holomorphic line bundles of the 
same degree{\rm )}. Moreover,  in this case $M$ admits a locally-symmetric 
CSC K\"ahler metric in each K\"ahler class.
\end{thm} 

Theorem~\ref{main} has already been proved 
under some 
additional assumptions. Indeed, Calabi \cite{cal-two}  described all 
extremal K\"ahler metrics on rational K\"ahler surfaces and in particular confirmed 
Theorem~\ref{main} when the genus ${\rm {\bf g}}=0$;  Burns and 
deBartolomeis \cite{BdB} proved Theorem~\ref{main} for K\"ahler 
classes $\Omega$ such that $c_1(M) \cdot \Omega =0$, Fujiki 
\cite{fujiki} confirmed it for ruled surfaces of genus ${\rm {\bf 
g}}=1$ while LeBrun \cite{LeB}  gave a proof when ${\rm {\bf g}} \ge 
2$ and the  K\"ahler class satisfies $c_1(M) \cdot \Omega \le 0$. 
Finally, in \cite{christina-three} it is shown that 
Theorem~\ref{main}  holds if the CSC K\"ahler metric is an Yamabe 
minimizer. 

To the best of our knowledge, no complete proof of Theorem~\ref{main} 
is readily available for the remaining case ({\bf g} $\ge$ 2 and 
$c_1(M)\cdot \Omega >0$). The purpose of this note is to fill this 
gap. Note that our arguments work only for ${\rm {\bf g}} \ge 2$, but do not require any assumption on 
the K\"ahler class.

\vspace{0.2cm}
As a standard corollary of Theorem~\ref{main} we also  derive 
(compare with \cite{LeB} and \cite[Proposition~4]{christina-one}):
\begin{cor} \label{standard} Let $(M, J) = P(E)$ be a geometrically 
ruled complex surface which admits an extremal K\"ahler metric. Then 
$E$ is either polystable {\rm (}in which case the scalar curvature is 
constant{\rm )} or else is the direct sum $L_{1} \oplus L_{2}$ of a pair of 
holomorphic line bundles of different degrees {\rm (}in which case the 
scalar curvature is not constant{\rm )}.
\end{cor}

\section{Proof of Theorem~\ref{main}} As we have already noted in 
Remark 1~(1), the `if' part of the theorem is well-known. So we deal 
with the `only if' part.  Let $(M,J) = P(E)$, where $\pi : E \to 
\Sigma_{\rm {\bf g}}$ is a 
holomorphic rank 2 vector bundle over a compact curve $\Sigma_{\rm 
{\bf g}}$ of genus {\bf g}. In the cases  {\bf g} $=0,1$, the theorem follows from \cite{cal-two} and \cite{fujiki}, respectively.  We will thus assume that {\bf g} $\ge 2$, and that $E$ is not 
polystable. Suppose that $(M,J)$ admits a CSC K\"ahler metric $g$ in 
a K\"ahler class $\Omega$. We want to arrive at a contradiction. 

Our first step is to verify the classical obstructions coming from 
the Matsushima--Lichnerowicz theorem~\cite{lichne,matsushima} and the vanishing of the Futaki 
invariant~\cite{futaki}; this is essentially done in \cite[Proposition~4]{christina-one}, but we outline the argument for seek of completeness.
\begin{lemma}\label{classic} If $(M,J)$ is as above, then $E$ is 
undecomposable and  ${\rm Aut}(M,J)$ is discrete.
In particular,  the only endomorphism of E is the scalar multiplication, 
i.e.,  E is {\rm simple}.

\end{lemma}
\begin{proof} We show that the Lie algebra $\mathfrak{h}$ of 
Aut$(M,J)$ of all  holomorphic vector fields on $(M,J)$ is trivial. 
This would imply that $E$ is undecomposable because otherwise 
multiplication with elements of $\C^{\times}$ of one of the split factors would induce a non-trivial 
holomorphic ${\mathbb C}^{\times}$-action on $(M,J)$. 

Suppose for 
a contradiction that  ${\mathfrak h} \neq \{ 0 \}$.  

Note that, since the 
Euler characteristic of $M$  is $e(M)=4(1-{\rm {\bf g}}) < 0$,  any 
holomorphic 
vector field has zeroes.  Therefore, by the Matsushima--Lichnerowicz theorem~\cite{lichne,matsushima}, 
${\mathfrak h}$ is the complexification of the real Lie algebra 
of Killing vector fields of the CSC K\"ahler metric $g$.  Thus, one can pick an $S^1$ isometric 
action which generates a ${\mathbb C}^{\times}$ holomorphic action on 
$(M,J)$, corresponding to an element $\Theta \in {\mathfrak h}$. 

The 
Lie algebra $\mathfrak h$ of $(M,J)$  is described by Maruyama 
\cite{mar}. In our case,  $\mathfrak h \cong H^0(\Sigma_{\rm {\bf g}}, {\rm sl}(E))$ so that 
for an element  $\Theta$ to generate a $\C^{\times}$-action it must be a (fiber-wise)
diagonalizable holomorphic section of the bundle ${\rm sl}(E)$. The coefficients of the characteristic polynomial of $\Theta$ are holomorphic functions on $\Sigma_{\rm {\bf g}}$, and therefore constants. It then follows that $\Theta$ defines a direct sum decomposition $E= L_1 \oplus L_{2}$, 
where the holomorphic sub-bundles $L_1$ and $L_2$ correspond to the 
eigenspaces of $\Theta$ at each fiber. Tensoring $E$ with $L_{1}^{-1}$ does not change the 
biholomorphic type of $M$, so without loss we can assume that $E = 
{\mathcal O} \oplus L$ and that $\Theta$ is the holomorphic vector field generated by $\C^{\times}$ multiplications of $L$. 

Since $E$ is 
not polystable by assumption,  ${\rm deg}(L)$ must be non-zero, and  using $P(\cO \oplus L) \cong P(\cO\oplus L^{-1})$, we can assume without loss that ${\rm deg}(L) >0$.
In this case, $c_1(L) = [\omega_{h}/2\pi]$ for a Hodge K\"ahler class $[\omega_h]$ of a CSC K\"ahler metric $h$ on $\Sigma_{\rm {\bf g}}$. Following~\cite{cal-two,Le-Sim,hwang-singer,ACGCT}, one can then construct a family of explicit K\"ahler metrics on $M$, which exhaust the K\"ahler cone of $M=P(\cO \oplus L)$. Picking such a metric in $\Omega$,  one easily computes the Futaki invariant $\frak{F}_{\Omega}(\Theta)$ of the holomorphic field $\Theta$ (see \cite[Theorem~3]{Le-Sim}, \cite[p.181] {christina-one} or \cite[\S 2.3]{ACGCT} for a precise formula). It turns out 
that $\frak{F}_{\Omega}(\Theta)$ is non-zero,   which  contradicts the existence of CSC K\"ahler metric on $M$~\cite{futaki}. \end{proof}

Recall that, by Narasimhan--Ramanan approximation 
theorem~\cite[Proposition~2.6]{NR}, any holomorphic vector bundle $E$ over 
$\Sigma_{\rm {\bf 
g}}$ (${\rm {\bf g}} \ge 2$) can be included in an analytic family of 
vector bundles over $\Sigma_{\rm {\bf g}}$,  in which stable vector 
bundles form an open dense subset. Such a family defines {\it stable 
deformations} of $E$. We can now apply Theorem~3 from 
\cite{fujiki}:
\begin{thm}\label{deformation}{\rm \cite{fujiki}} Let $(M,J)=P(E)$ be a 
geometrically ruled manifold over a compact curve $\Sigma_{\rm {\bf 
g}}$ with  ${\rm {\bf g}} \ge 2$,  such that $E$  is simple and 
non-polystable. Suppose 
that $M$ admits a CSC K\"ahler metric.  Then there exists a small 
stable 
deformation  $E'$ of $E$, such that the ruled manifold $M' = P(E')$ 
admits an 
extremal K\"ahler metric different from the locally-symmetric one.  
\end{thm}

But Theorem~\ref{deformation} contradicts the uniqueness established 
by Theorem~\ref{uniqueness}, which completes the proof of our
Theorem~\ref{main}. 

\section{Proof of Corollary~\ref{standard}} 
 By Theorem~\ref{main}, we may assume that the scalar curvature $s_g$ 
of our extremal metric $g$ is not constant. Moreover, the case when 
$E$ is a vector bundle over ${\mathbb C}P^1$ (i.e. $M$ is a rational 
ruled surface) is completely solved by Calabi~\cite{cal-two}.  We thus 
may also assume that $E$ is  a vector bundle over a Riemann surface 
$\Sigma_{\rm {\bf g}}$ of  genus ${\rm {\bf g}} \ge 1$.

The isometry group of $g$ is non-trivial, since $J {\rm grad}_{g} 
s_g$ is a Killing field. Since this Killing field has zeroes, there 
exists
a non-trivial $S^1$-action by isometric biholomorphisms, and with 
fixed points. 
Since ${\rm {\bf g}} \ge 1$ (and therefore the base $\Sigma_{\rm {\bf 
g}}$ admits no non-trivial holomorphic fields with zeroes),  this 
$S^1$-action preserves the ${\mathbb C}P^1$-fibers, and thus comes 
from a diagonalizable holomorphic section of ${\rm sl}(E)$ (see the 
proof of Lemma~\ref{classic}), which defines the  desired direct sum 
decomposition of $E$.

\section{Concluding remarks}

(1) The arguments used in the proof of Theorem~\ref{main} imply 
non-existence of CSC K\"ahler metrics on the projectivization of an {\it undecomposable} unstable  vector bundle of {\it arbitrary rank}, over a  compact complex curve $\Sigma_{\bf g}$ of genus ${\rm {\bf g}} \ge 
2$.  Note, however, that in this more general setting we don't have 
an analogue of the second part of Lemma~\ref{classic}, nor can we directly appeal to 
Narasimhan--Ramanan approximation theorem~\cite{NR} when ${\rm {\bf g}}=1$. On 
the other hand, as we have already explained in the introduction, the 
results of \cite{RT} and \cite{CT,Do-one,mab-one} show non-existence of 
CSC K\"ahler metrics in the rational K\"ahler classes of 
$P(E)$, for $E$ non-semistable and of 
arbitrary rank, and ${\rm {\bf g}} \ge 1$. 

(2) The existence issue for extremal K\"ahler metrics of non-constant 
scalar curvature on geometrically ruled complex manifolds seems to be a more 
delicate problem. In the cases when ${\rm {\bf g}} =0,1$ it is known that any K\"ahler class admits an extremal K\"ahler metric ~\cite{cal-two,hwang,christina-two};  however,  when  ${\rm {\bf g}}>1$, by combining Theorem~\ref{uniqueness} with an 
observation from \cite{christina-one}, one 
obtains examples of ruled surfaces which exhibit both existence and 
non-existence for different K\"ahler classes. This phenomenon is related in \cite{szekelyhidi} to an appropriate equivariant $K$-polystability notion.

\vspace{0.3cm}
\noindent
{\sc Acknowledgements.} We would like to thank C.~LeBrun for encouraging us to write 
this note, J.~Ross and R.~Thomas  for explaining to us their work and sharing with us their expertise on various stability notions, and M.~Singer and the referee for valuable remarks.

\end{document}